\documentclass[a4paper,10pt]{article}
\usepackage{amsmath}
\usepackage{amssymb}
\usepackage{epsf}
\usepackage{epsfig}
\usepackage{multirow,longtable}
\usepackage[latin1]{inputenc}
\usepackage{amsfonts}
\usepackage{latexsym}
\usepackage{pstricks}
\usepackage{graphicx}

\usepackage[algo2e,lined,algoruled,titlenumbered,algosection]{algorithm2e}
\usepackage[T1]{fontenc}       

\newcommand{\dem}{\par \noindent{\bf Proof:} \par}
\newcommand{\fin}{\hfill $\square$  \par \bigskip}
\newtheorem{teorema}{\bf Theorem}[section]
\newtheorem{lema}{\bf Lemma}[section]
\newtheorem{definicion}{\bf Definition}[section]
\newtheorem{proposicion}{\bf Proposition}[section]

\newtheorem{observacion}{\bf Remark}[section]

\begin{document}

\title{On the simultaneous location of a service facility and a rapid transit line}
\author{I. Espejo \and  A.M. Rodríguez-Chía   }
\date{Dpto. Estadística e Investigación Operativa. Universidad de Cádiz }
\maketitle \thispagestyle{empty}
\begin{abstract}
In this paper we provide a finite set of candidates to be one of the endpoints of
an optimal solution  for the problem of locating simultaneously a service facility and a rapid transit line.
\end{abstract}
\section{Introduction}

The problem of locating simultaneously a service facility and a rapid transit line was introduced by  
Espejo and Rodríguez-Chía \cite{ERCH11}. In that paper, the authors present  a set of candidates to be  
one of the endpoints of an optimal transit line as well as  a solution procedure for this problem. Later, 
Díaz-Bañez et al.  \cite{DB11} provided an alternative description of the  endpoints of an optimal solution 
for this problem because the previous description does not work in general. In this paper, we give a finer 
description of this set, actually we provide a set of candidates of finite cardinality to be one of the endnodes of the optimal segment.

A complete motivation and justification of above mentioned model is given in \cite{ERCH11}. Therefore, we proceed directly with the formulation of the problem.
Formally, consider $(\mathbb{R}^2,\,||\cdot ||_1)$ and let
$A=\{a_1,\ldots,a_M\}$ be the set of demand points,
$\omega=\{\omega_1,\ldots,\omega_M\}\subseteq \mathbb{R}^+$ be the
set representing the intensity of the demand in the elements of the
set $A$ and $x$ be the service facility to be located. The rapid transit line, represented by a segment of given
length, $\ell$, is defined by two extreme points, $e$ and $s$, which
represent the entrance (access) and the exit, respectively. {The
distance from $e$ to $s$ through the rapid transit line is defined
as} $\ell/k$, where $\ell$ is the length of the rapid transit line
and $k$ is a constant greater than or equal to 1. (Note that it
means that the speed within the rapid transit line is greater than
or equal to in the plane).
Obviously, the closest extreme point of the segment to the service
facility should be the one representing the exit, that is,
$||x-s||_1\leq ||x-e||_1.$
In \cite{ERCH11} is proved that an  optimal location for service facility is placed on the exit extreme point of the optimal
segment for the problem of locating simultaneously a service facility and a rapid transit line.
Therefore, the model can be reformulated as the problem of locating a rapid transit line (segment) such that one its endpoints is the service facility.
Thus, the distance between a demand point and the service facility is
given by the function
$$d_{e,x}(a_i,x)=\min\{||x-a_i||_1,||e-a_i||_1+\ell/k\}.$$
Therefore, the  problem can be formulated as follows:
\begin{eqnarray}\label{mod_general}
&\min &f(e,x)=\sum_{i=1}^Mw_i d_{e,x}(a_i,x)\\ \nonumber
&s.t.&||x-e||_2=\ell.
\end{eqnarray}

Observe that for $\ell=0$, the classical (weighted) rectilinear
median problem in the plane is obtained. Let $W$ be the solution set
of this problem for the case $\ell=0$  (usually called median, $me$). Let $I$ be the set of intersection points (intersection points of the fundamental lines,
in this case, horizontal and vertical lines at demand points, see
\cite{DurierM:85}) and $I_W:=I\cap W$. 
An interesting result obtained in \cite{ERCH11}  states  that the endpoints of  an optimal segment are in opposite quadrants with respect to some $t\in I_W$, see \cite{ERCH11} for further details.

In general, the four quadrants defined by a point $q\in\mathbb{R}^2$
will be denoted by $C^i_q$, $i=1,\ldots,4$ (to simplify the
notation, the quadrants defined by the origin will be denoted by
$C^i$, $i=1,\ldots,4$).
In what follows, we study the case $|I_W|=1$ (notice that in this
case $|W|=1$, that is, the median is unique), the remaining cases
will be analyzed analogously. Without loss of generality and for a better
understanding, it is assumed that the median point is the origin,
i.e., $I_W=me={(0,0)}$ and in addition, we solve the problem where $e$ is
located at the first quadrant defined by the origin, $C^1$ (the
analysis of the remaining cases can be analogously addressed). Hence, $x$
will be in $C^3$.
In general, the circumference centered at point $x$ of radius $\ell$ will be denoted
by $S_x$ (to simplify the notation, the circumference centered at
the origin will be denoted by $S$). Let $\overline{S}$ denote the
circle centered at the origin and radius $\ell$,
$\overline{S}\,^i:=\overline{S}\cap C^i$, $S^i_x:=S_x\cap C^i$,
$C^i(A):=C^i\cap A$, {for ${i=1,2,3,4}$}.

\begin{figure}[htb]
\begin{center}
\includegraphics[width=13cm]{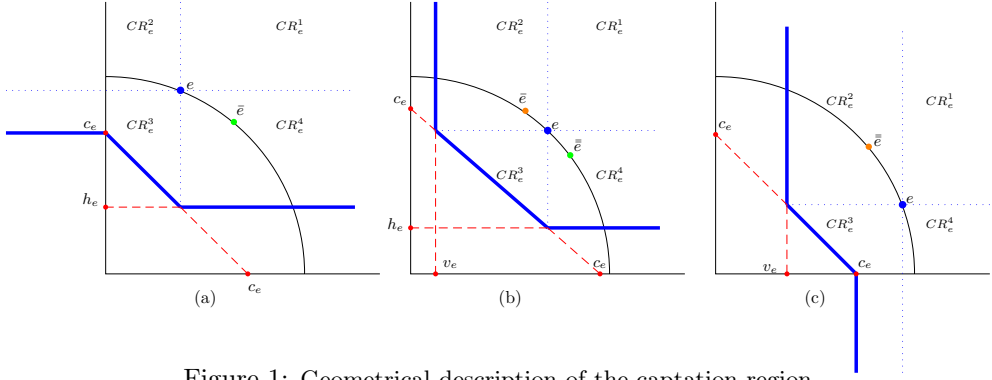}
\end{center}\vspace*{-1.1cm}
\caption{\small Geometrical description of the captation region}
\label{rc_3}
\end{figure}

In order to analyze Problem \eqref{mod_general}, we recall the definition of captation region:
\begin{definicion}\label{def_rc}
{For  a given segment with endpoints $e$ and $x$, the set $CR_{e,x}:=\{a=(a_1,a_2)\in A: ||x-a||_1\geq
||e-a||_1+\ell/k \}$ is called captation region associated with $e$.
The boundary of $CR_{e,x}$ {is defined as} $\partial CR_{e,x}=\{a\in CR_{e,x}:
||x-a||_1= ||e-a||_1+\ell/k\}$
and the relative interior of $CR_e$ as $\mbox{ri}(CR_{e,x})= CR_{e,x}\setminus \partial CR_{e,x}.$ }\end{definicion} 
Note that  $a_i\in CR_e$ if and only if $d_e(x,a_i)=\min{\{||x-a_i||_1,
||e-a_i||_1+\ell/k\}}=||e-a_i||_1+\ell/k.$ Therefore, all demand
points that belong to the captation region associated with $e\in
S^1$ will use the rapid transit line to achieve the service
facility. The captation region a\-sso\-cia\-ted with $e$
and $x$ is given by (see \cite{ERCH11} for further details):

\noindent - if $x_1\leq e_1\leq \bar{e}_1$,
$$CR_{e,x}=\{a=(a_1,a_2)\in A:
a_1\geq x_1,\; a_2\geq h^+_{e,x}, \; a_1+a_2\geq c_{e,x} \}\cup
\{a=(a_1,a_2)\in A: a_1<x_1,\; a_2\geq h^-_{e,x} \}, $$
- if $\bar{e}_1<e_1<\bar{\bar{e}}_1$,
$$CR_{e,x}=\{a=(a_1,a_2)\in
A: a_1\geq v^+_{e,x},\; a_2\geq h^+_{e,x}, \; a_1+a_2\geq c_{e,x}
\},
$$
- if $\bar{\bar{e}}_1\leq e_1\leq x_1+\ell$,
$$CR_{e,x}=\{a=(a_1,a_2)\in
A: a_2\geq x_2,\; a_1\geq v^+_{e,x}, \; a_1+a_2\geq c_{e,x} \}\cup
\{a=(a_1,a_2) \in A: a_2<x_2,\; a_1\geq v^-_{e,x} \},$$
where $h^-_{e,x}:=\frac{e_1+e_2+\ell/k}{2}+\frac{x_2-x_1}{2}$,
$h^+_{e,x}:=\frac{-e_1+e_2+\ell/k}{2}+\frac{x_1+x_2}{2}$,
$v^-_{e,x}:=\frac{e_1+e_2+\ell/k}{2}+\frac{x_1-x_2}{2}$,
$v^+_{e,x}:=\frac{e_1-e_2+\ell/k}{2}+\frac{x_1+x_2}{2}$,
$c_{e,x}:=\frac{e_1+e_2+\ell/k}{2}+\frac{x_1+x_2}{2}$ and
$\bar{e}=(\bar{e}_1,\bar{e}_2)$,
$\bar{\bar{e}}=(\bar{\bar{e}}_1,\bar{\bar{e}}_2) \in S_x$ such that
$(\bar{e}_2-x_2)-(\bar{e}_1-x_1)=\frac{\ell}{k}$
($\bar{e}_2-x_2>\bar{e}_1-x_1$) and
$(\bar{\bar{e}}_1-x_1)-(\bar{\bar{e}}_2-x_2)=\frac{\ell}{k}$
($\bar{\bar{e}}_1-x_1>\bar{\bar{e}}_2-x_2$). 
%
%
The boundary of $CR_{e,x}$ is defined as $\partial CR_{e,x}=\{a\in
CR_{e,x}: ||a-x||_1=||a-e||_1+\ell/k\}$. A graphical
representation of these captation regions for the case $(x_1,x_2)=(0,0)$ is given by Figure
\ref{rc_3}. 

\section{Dominant solutions}

In this section, we characterize a   finite set of candidates
to be one of the endpoints of an optimal solution  of Problem (\ref{mod_general}).

\begin{figure}[h]
\begin{center}
\includegraphics[width=4.7cm]{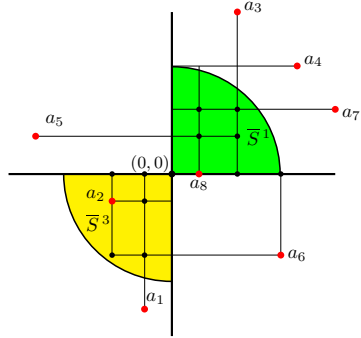}
\vspace{-0.5cm}\caption{\small Possible locations of $e$ and $x$.
\label{Gceldas1}} \end{center}
\end{figure}

\begin{teorema}
\label{t21}
Let $e^*$ and $x^*$ be the endpoints of an optimal
segment of Problem (\ref{mod_general}) with $e^*$ restricted to the
first quadrant. Then,
 $x^*$ and $e^*$, satisfy one of the following conditions:
\begin{itemize}
\item[i)] Either $x^*\in \bar{S}^3\cap I$ or $e^*\in
\bar{S}^1\cap I$ (Figure \ref{Gceldas1} shows an example of the
possible locations of the endpoints of an optimal segment).
\item[ii)] $x^*_k=a_{jk}$ and $e^*_{k'}=a_{j'k'}$ for some $j,j' \in \{1,\ldots, M\}$ and
$k,k' (\ne k) \in \{1,2\}$, such that, the angle of the segment with the positive direction of the $x$-edge, $\theta$, satisfies that
$\tan (\theta)= -\frac{w^b_{e^*,x^*}}{w^a_{e^*,x^*}}$ with
$w^b_{e^*,x^*}>0$ (and consequently $w^a_{e^*,x^*}<0$), where
$w^a_{e^*,x^*}$ and $w^b_{e^*,x^*}$ are defined in the proof of Lemma \ref{prop42_2}.
\end{itemize}
\end{teorema}

\begin{observacion}
The description given by \cite{DB11} states that either one of the endpoints of an optimal segment is on an intersection point or both endpoints are placed on fundamental lines at two demand points, one vertical and the other horizontal. Observe that we have with this description a infinite many  number of candidates to be one of the endpoints of an optimal segment. However, our description is finer than the one in \cite{DB11} because we have reduced this set of candidates to  a finite set. 
\end{observacion}

In order to prove this result, we will show that for any segment of length
$\ell$ defined by $x\not\in \bar{S}^3\cap I$ and $e\not \in
\bar{S}^1\cap I$, {such that, $x$ and $e$ do not satisfy condition ii) 
of Theorem \ref{t21},} we can find another segment of extreme points $x'$
and $e'$ such that $f(x',e')\leq f(x,e)$.
Observe that if $x\not\in \bar{S}^3\cap I$ ($e\not \in \bar{S}^1\cap
I$), then $x_1\neq a_{i1}$ or $x_2\neq a_{i2}$, $\forall a_i\in A$
($e_1\neq a_{i1}$ or $e_2\neq a_{i2}$, $\forall a_i\in A$). Firstly,
the case $x_1\neq a_{i1}$ and $e_1\neq a_{i1}$, $\forall a_i\in A$
will be considered (Subsection \ref{sub1}); secondly, $x_1\neq
a_{i1}$ and $e_2\neq a_{i2}$, $\forall a_i\in A$ (Subsection
\ref{sub2}). The remaining cases can be analogously studied. 

\subsection{Case $x_1\neq a_{i1}$ and $e_1\neq
a_{i1}$, $\forall a_i\in A$.}\label{sub1}

\begin{figure}[h]
\centerline{\parbox{210pt}{\begin{center}
\includegraphics[width=5.8cm]{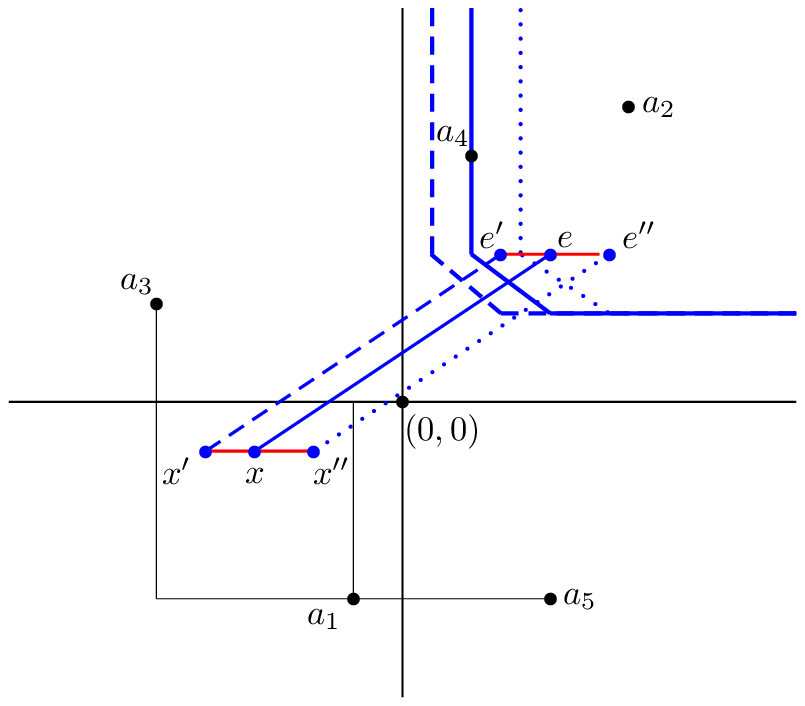}
\vspace{-0.1cm} \caption{\small $e_1\neq a_{i1}$, $\forall a_i\in
A$. \label{seg_may_2}}
\end{center}
}
\parbox{210pt}{\begin{center}
\includegraphics[width=6cm]{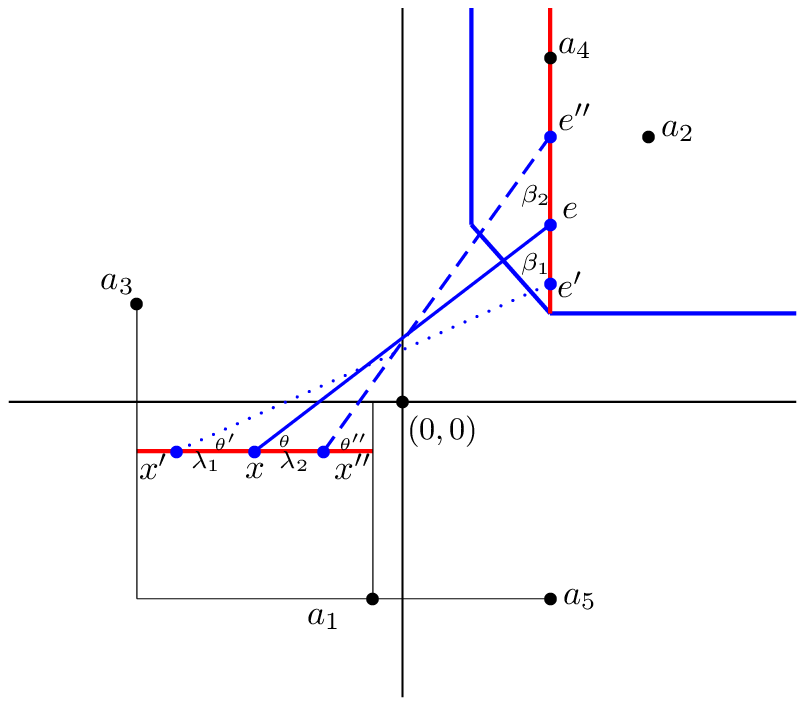}
\vspace{-0.7cm}\caption{\small $e_2\neq a_{i2}$, $\forall a_i\in
A$.\label{seg_in_cuad_2}}
\end{center}
}}
\end{figure}

Let $x'$, $e'$, $x''$ and $e''$ be such that (see Figure
\ref{seg_may_2})
\begin{eqnarray}\label{desplaza1_2}
x'=x+(-\lambda,0);\quad e'=e+(-\lambda,0);&& x''=x+(\lambda,0);\quad
e''=e+(\lambda,0),
\end{eqnarray}
where $\lambda$ is a small enough positive value satisfying that:
\begin{itemize}
\item[i)] $CR_{e,x}\backslash \partial CR_{e,x}=CR_{e',x'}
\backslash \partial CR_{e,x}=CR_{e'',x''}
\backslash \partial CR_{e,x}$ and $C^i_x(A)
=C^i_{x'}(A)
=C^i_{x''}(A)
$,
$i=1,\ldots,4$.
\item[ii)] $a_{i1}\neq e'_1$ and $a_{i1}\neq e''_1$, for any $a_{i}\in
A$.
\end{itemize}
\begin{teorema}\label{teo1_2}
Let $e$, $x$ be the extreme points of a segment of length $\ell$
such that $x_{1}\neq a_{i1}$ and $e_1\neq a_{i1}$, $\forall a_i\in
A$. Then, $f(e',x')\leq f(e,x)$ or $f(e'',x'')\leq f(e,x)$, where
$e'$, $x'$, $e''$, $x''$ were given by (\ref{desplaza1_2}).
\end{teorema}
\dem
Since $||x-a_i||_1\leq ||a_i-e||_1+\ell/k$, for all $a_i\in
C^3_x(A)$, the objective function can be expressed as:
\begin{eqnarray*}
f(e,x)\!\!&=&\!\!\!\!\!\!\!\!\!\!\!\!\!\!\sum_{a_i\not\in
CR_{e,x}}\!\!\!\!\!\!\!\!\omega_i||x-a_i||_1+
\!\!\!\!\!\!\!\!\!\!\sum_{a_i\in CR_{e,x}}\!\!\!\!\!\!\!\!\omega_i
(||e-a_i||_1+\ell/k)=\!\!\!\!\!\!\!\!\!\!\!\!\sum_{a_i\in
(C^2_{x}(A)\backslash CR_{e,x})\cup
C^3_{x}(A)}\!\!\!\!\!\!\!\!\!\!\!\!\!\omega_i
||x-a_i||_1+\!\!\!\!\!\!\!\!\!\!\!\!\!\!\!\!\!\!\!\!\!\sum_{a_i\in
(C^1_{x}(A)\cup C^4_{x}(A))\backslash
CR_{e,x}}\!\!\!\!\!\!\!\!\!\!\!\!\!\!\!\!\!\!\omega_i
||x-a_i||_1\\
\!\!&+&\!\!\!\!\!\!\!\!\!\!\sum_{a_i\in
CR_{e,x}\backslash \partial CR_{e,x}}\!\!\!\!\!\omega_i
(||e-a_i||_1+\ell/k)+\!\!\!\!\!\!\sum_{a_i\in
\partial CR_{e,x}}\!\!\!\!\!\omega_i (||e-a_i||_1+\ell/k).
\end{eqnarray*}

From (\ref{desplaza1_2}), the definition of $\lambda$ and $e_1\neq
a_{i1}$, $\forall a_i\in A$, we obtain the following:
\begin{enumerate}
\item[(1)]$\forall a_i\in C^2_x(A)\backslash \mbox{ri}( CR_{e,x})\cup C^3_x(A)$:
$||x'-a_{i}||_1=||x-a_i||_1-\lambda$ and
$||x''-a_i||_1=||x-a_i||_1+\lambda$.
\item[(2)]$\forall a_i\in (C^1_x(A)\cup C^4_x(A))\backslash CR_{e,x})$: $||x'-a_{i}||_1=||x-a_i||_1+\lambda$ and
$||x''-a_i||_1=||x-a_i||_1-\lambda$.
\item[(3)]$\forall a_i\in$ $CR_{e,x}\backslash \partial CR_{e,x}$:
\begin{itemize}
\item[(a)]$a_{i1}<e_1$: $||e'-a_i||_1=||e-a_i||_1-\lambda$ and
$||e''-a_i||_1=||e-a_i||_1+\lambda$.
\item[(b)] $a_{i1}>e_1$:
$||e'-a_i||_1=||e-a_i||_1+\lambda$ and
$||e''-a_i||_1=||e-a_i||_1-\lambda$.
\end{itemize}
\item[(4)]$\forall a_i\in \partial CR_{e,x}\setminus C_{x}^2(A)$:
\begin{itemize}
\item[(a)] $a_{i1}>e_1$:
$||e'-a_i||_1+\ell/k=||e-a_i||_1+\lambda+\ell/k=||x-a_i||_1+\lambda=||x'-a_i||_1$
and
$||e''-a_i||_1+\ell/k=||e-a_i||_1-\lambda+\ell/k=||x-a_i||_1-\lambda=||x''-a_i||_1$.
\item[(b)]$a_{i1}<e_1$: $||e'-a_i||_1+\ell/k=||e-a_i||_1-\lambda+\ell/k=||x-a_i||_1-\lambda < ||x-a_i||_1+\lambda=||x'-a_i||_1$
and
$||e''-a_i||_1+\ell/k=||e-a_i||_1+\lambda+\ell/k=||x-a_i||_1+\lambda
>||x-a_i||_1-\lambda=||x''-a_i||_1$.
\end{itemize}
\end{enumerate}
Hence, we get that
\begin{eqnarray*}
f(x',e')-f(x,e)&=&\left(\omega^1_{e,x}-\omega^2_{e,x}-\omega^3_{e,x}+\omega^4_{e,x}+\omega^{+}_{e,x}-\omega^{-}_{e,x}
+\omega^{\delta +}_{e,x} - \omega^{\delta-}_{e,x}\right) \lambda =
\Delta',\\
f(x'',e'')-f(x,e)&=&\left(-\omega^1_{e,x}+\omega^2_{e,x}+\omega^3_{e,x}-\omega^4_{e,x}-\omega^{+}_{e,x}+\omega^{-}_{e,x}
-\omega^{\delta+}_{e,x}-\omega^{\delta-}_{e,x}\right)\lambda=\Delta'',
\end{eqnarray*}
where
$$\begin{array}{lclllclllcll}
\omega^1_{e,x}&=&\!\!\!\!\!\!\!\!\!\!\!\!\!\!\!\!\!\!\displaystyle
\sum_{a_j\in C^1_x(A)\backslash
CR_{e,x}}\!\!\!\!\!\!\!\!\!\!\!\!\!\!\omega_j;&
\omega^{2}_{e,x}&=&\!\!\!\!\!\!\!\!\!\displaystyle\sum_{a_j\in
C^2_x(A)\backslash \mbox{ri}(CR_{e,x})}\!\!\!\omega_j;&
\omega^{3}_{e,x}&=&\!\!\!\!\!\!\!\displaystyle \sum_{a_j\in
C^3_x(A)}\!\!\omega_j;&
\omega^{4}_{e,x}&=&\!\!\!\!\!\!\!\!\displaystyle \sum_{a_j\in
C^4_x(A)\backslash CR_{e,x}}\!\!\!\!\omega_j;\\
\omega^{+}_{e,x}&=&\!\!\!\!\!\!\!\!\!\!\!\!\!\!\!\!\!\!\displaystyle
\sum_{a_j\in CR_{e,x}\backslash \partial CR_{e,x}:
a_{j1}>e_1}\!\!\!\!\!\!\!\!\!\!\!\!\!\!\!\!\!\omega_j;&
\omega^{-}_{e,x}&=&\!\!\!\!\!\!\!\!\!\!\!\!\!\!\!\!\!\!\displaystyle
\sum_{a_j\in CR_{e,x}\backslash \partial CR_{e,x}:
a_{j1}<e_1}\!\!\!\!\!\!\!\!\!\!\!\!\!\!\!\!\!\omega_j;&
\omega^{\delta
+}_{e,x}&=&\!\!\!\!\!\!\!\!\!\!\!\!\!\!\!\!\!\!\displaystyle
\sum_{a_j\in \partial CR_{e,x}\setminus C_x^2(A):
a_{j1}>e_1}\!\!\!\!\!\!\!\!\!\!\!\!\!\!\!\!\!\omega_j;&
\omega^{\delta
-}_{e,x}&=&\!\!\!\!\!\!\!\!\!\!\!\!\!\!\!\!\!\!\displaystyle
\sum_{a_j\in \partial CR_{e,x}\setminus C_x^2(A):
a_{j1}<e_1}\!\!\!\!\!\!\!\!\!\!\!\!\!\!\!\!\!\omega_j.
\end{array}$$
Since $\Delta'=-\Delta''- 2\lambda\omega^{\delta-}_{e,x}$, we
conclude that:
\begin{itemize}
\item[-]If $\Delta''\leq 0$, then $f(x'',e'')\leq f(x,e)$.
\item[-]Otherwise, $\Delta'\leq 0$ and $f(x',e')\leq f(x,e)$,
and the result follows.\fin
\end{itemize}

\subsection{ Case $x_1\neq a_{i1}$  and $e_2\neq a_{i2}$, $\forall
a_i\in A$.}\label{sub2}

Let $x'$, $e'$, $x''$ and $e''$ be such that (see Figure
\ref{seg_in_cuad_2})
\begin{eqnarray}\label{desplaza_2}
x'=x+(-\lambda_1,0);\quad e'=e+(0,-\beta_1);&& x''=x+(\lambda_2,0);
\quad e''=e+(0,\beta_2),
\end{eqnarray} where
$\lambda_1, \lambda_2, \beta_1, \beta_2\in \mathbb{R}^+ $ are small
enough values satisfying that
\begin{itemize}
\item[i)]$CR_{e,x}\backslash \partial CR_{e,x}=CR_{e',x'}
\backslash \partial CR_{e,x}=CR_{e'',x''}\backslash \partial
CR_{e,x}$ and $C^i_x(A)
=C^i_{x'}(A)
=C^i_{x''}(A)
$,
$i=1,\ldots,4$.
\item[ii)] $e'_2\neq a_{i2}$ and $e''_2\neq a_{i2}$, for any $a_i\in
A$.
\end{itemize}
In order to prove that either $f(e',x')\leq f(e,x)$ or
$f(e'',x'')\leq f(e,x)$, we will use similar arguments as in the
previous case, before that we give some technical results. 

\begin{proposicion}\label{observt}
Let $\theta$, $\theta'$ and $\theta''$ be the angles defined by the
segment of extreme points $e$ and $x$, $e'$ and $x'$, $e''$ and
$x''$, and the horizontal line, respectively (see Figure
\ref{seg_in_cuad_2}), where $e'$, $x'$, $e''$ and $x''$ were defined
in (\ref{desplaza_2}). Then
\begin{eqnarray*}
\frac{\beta_1}{\lambda_1}=\displaystyle\frac{2\ell\sin\left(\frac{\theta-\theta'}{2}\right)
\cos\left(\frac{\theta+\theta'}{2}\right)}
{2\ell\sin\left(\frac{\theta'+\theta}{2}\right)
\sin\left(\frac{\theta-\theta'}{2}\right)}
=\displaystyle\frac{1}{\tan(\frac{\theta+\theta'}{2})},&&
\frac{\beta_2}{\lambda_2}=\displaystyle\frac{2\ell\sin\left(\frac{\theta''-\theta}{2}\right)
\cos\left(\frac{\theta''+\theta}{2}\right)}
{2\ell\sin\left(\frac{\theta+\theta''}{2}\right)
\sin\left(\frac{\theta''-\theta}{2}\right)}
=\displaystyle\frac{1}{\tan(\frac{\theta''+\theta}{2})}.
\end{eqnarray*}
Moreover,
\begin{itemize}
\item[i)]If $\lambda_1=\lambda_2$ ($\beta_1=\beta_2$) then
$\beta_2<\beta_1$ ($\lambda_2>\lambda_1$).
\item[ii)]If $\theta'<\theta<\theta''\leq \pi /4$ ($\pi /4\leq
\theta'<\theta<\theta''$), then $\lambda_1<\beta_1$ and
$\lambda_2<\beta_2$ ($\lambda_1>\beta_1$ and $\lambda_2>\beta_2$).
\end{itemize}
\end{proposicion}
\dem

It is consequence of $ \lambda_1=\ell(\cos{\theta'}-\cos{\theta})$,
$\beta_1=\ell(\sin{\theta}-\sin{\theta'})$,
$\lambda_2=\ell(\cos{\theta}-\cos{\theta''})$,
$\beta_2=\ell(\sin{\theta''}-\sin{\theta})$, and that
$0<\theta'<\theta<\theta''\leq \frac{\pi}{2}$ and the tangent
function is monotonous in $(-\frac{\pi}{2}, \frac{\pi}{2})$. \fin


\begin{lema}\label{prop42_2}
Let $e'$, $x'$, $e''$, $x''$ be given by (\ref{desplaza_2}). Then,
$$f(x',e')-f(x,e)=\Delta'\mbox{ and } f(x'',e'')-f(x,e)=\Delta'',$$
where $\Delta'=\omega^a_{e,x}\lambda_1+\omega^b_{e,x}\beta_1+d'$ and
$\Delta''=-\omega^a_{e,x}\lambda_2-\omega^b_{e,x}\beta_2+d''$, with
$\omega^a_{e,x},\omega^b_{e,x},d',d''\in \mathbb{R}$, $d', d''\leq
0$.
\end{lema}
%
%

\dem
Let $e'$, $x'$, $e''$, $x''$ be given by (\ref{desplaza_2}).

\begin{enumerate}
\item If $a_i\in \partial CR_{e,x}\setminus C_x^2(A)$, with $a_{i2}\leq e_2$, then:
\begin{eqnarray*}
\!\!\!\!\!\!\!\!\!\!\!\!\!\!\!\!\!\!||e'-a_i||_1+\ell/k=||e-a_i||_1-\beta_1+\ell/k=||x-a_i||_1-\beta_1
< ||x-a_i||_1+\lambda_1=||x'-a_i||_1,&\mbox{i.e. } a_i\in
CR_{e',x'},&\\
\!\!\!\!\!\!\!\!\!\!\!\!\!\!\!\!\!\!||e''-a_i||_1+\ell/k=||e-a_i||_1+\beta_2+\ell/k=||x-a_i||_1+\beta_2
>||x-a_i||_1-\lambda_2=||x'-a_i||_1,&\mbox{ i.e. } a_i\not\in CR_{e'',x''}.&
\end{eqnarray*}

\item If $a_i\in \partial CR_{e,x}\setminus C_x^2(A)$, with $a_{i2}> e_2$, using similar
arguments and by Proposition  \ref{observt}, it can be proved that:
\begin{enumerate}
\item $\theta'<\theta<\theta''\leq \pi/4$, then $a_i\in CR_{e'',x''}$ and $a_i\not\in
CR_{e',x'}$.
\item $\pi/4\leq \theta'<\theta<\theta''$, then $a_i\not\in CR_{e'',x''}$ and $a_i\in
CR_{e',x'}$.
\item $\theta'<\theta=\pi/4<\theta''$, then $a_i\not\in CR_{e',x'}$ and $a_i\not\in
CR_{e'',x''}$.
\end{enumerate}
\end{enumerate}
Consider the following notation:
$$\begin{array}{lclllclllcll}
\omega^1_{e,x}&=&\!\!\!\!\!\!\!\!\!\!\displaystyle \sum_{a_j\in
C^1_x(A)\backslash CR_{e,x}}\!\!\!\!\!\!\!\omega_j;&
\omega^{2}_{e,x}&=&\!\!\!\!\!\!\!\!\!\displaystyle\sum_{a_j\in
C^2_x(A)\backslash \mbox{ri}( CR_{e,x})}\!\!\!\omega_j;&
\omega^{3}_{e,x}&=&\!\!\!\!\!\!\!\displaystyle \sum_{a_j\in
C^3_x(A)}\!\!\omega_j;&
\omega^{4}_{e,x}&=&\!\!\!\!\!\!\!\!\displaystyle \sum_{a_j\in
C^4_x(A)\backslash CR_{e,x}}\!\!\!\!\omega_j;
\\
\omega^{+}_{e,x}&=&\!\!\!\!\!\!\!\!\!\!\!\!\!\!\!\!\!\!\displaystyle
\sum_{a_j\in CR_{e,x}\backslash \partial CR_{e,x}:
a_{j2}>e_2}\!\!\!\!\!\!\!\!\!\!\!\!\!\!\!\!\!\omega_j;&
\omega^{-}_{e,x}&=&\!\!\!\!\!\!\!\!\!\!\!\!\!\!\!\!\!\!\displaystyle
\sum_{a_j\in CR_{e,x}\backslash \partial CR_{e,x}:
a_{j2}<e_2}\!\!\!\!\!\!\!\!\!\!\!\!\!\!\!\!\!\omega_j;&
\omega^{\delta
+}_{e,x}&=&\!\!\!\!\!\!\!\!\!\!\!\!\!\!\!\!\!\!\displaystyle
\sum_{a_j\in \partial CR_{e,x}\setminus C_x^2(A):
a_{j2}>e_2}\!\!\!\!\!\!\!\!\!\!\!\!\!\!\!\!\!\omega_j;&
\omega^{\delta
-}_{e,x}&=&\!\!\!\!\!\!\!\!\!\!\!\!\!\!\!\!\!\!\displaystyle
\sum_{a_j\in \partial CR_{e,x}\setminus C_x^2(A):
a_{j2}<e_2}\!\!\!\!\!\!\!\!\!\!\!\!\!\!\!\!\!\omega_j,
\end{array}$$

Using similar arguments as Theorem \ref{teo1_2}, we can obtain an
expression of $f(x',e')-f(x,e)=\Delta'$ and
$f(x'',e'')-f(x,e)=\Delta''$ as follows:
\begin{itemize}
\item[(a)]if $\theta'<\theta<\theta''\leq\frac{\pi}{4}$,
\begin{eqnarray*}
\Delta'&=&(\omega^1_{e,x}-\omega^2_{e,x}-\omega^3_{e,x}+\omega^4_{e,x})\lambda_1+
(\omega^{+}_{e,x}-\omega^{-}_{e,x})\beta_1
-\omega^{\delta-}_{e,x}\beta_1+\omega^{\delta+}_{e,x}\lambda_1\\
&=&(\omega^1_{e,x}-\omega^2_{e,x}-\omega^3_{e,x}+\omega^4_{e,x})\lambda_1+
(\omega^{+}_{e,x}-\omega^{-}_{e,x}+\omega^{\delta+}_{e,x})\beta_1
-\omega^{\delta-}_{e,x}\beta_1+\omega^{\delta+}_{e,x}(\lambda_1-\beta_1);\\
\Delta''&=&(-\omega^1_{e,x}+\omega^2_{e,x}+\omega^3_{e,x}-\omega^4_{e,x})\lambda_2+
(-\omega^{+}_{e,x}+\omega^{-}_{e,x})\beta_2
-\omega^{\delta-}_{e,x}\lambda_2-\omega^{\delta+}\beta_2\\
&=&(-\omega^1_{e,x}+\omega^2_{e,x}+\omega^3_{e,x}-\omega^4_{e,x})\lambda_2+
(-\omega^{+}_{e,x}+\omega^{-}_{e,x}-\omega^{\delta+}_{e,x})\beta_2
-\omega^{\delta-}_{e,x}\lambda_2.
\end{eqnarray*}
Therefore, $\Delta'$ and $\Delta''$ can be expressed as
$\Delta'=\omega^a_{e,x}\lambda_1+\omega^b_{e,x}\beta_1+d_1$ and
$\Delta''=-\omega^a_{e,x}\lambda_2-\omega^b_{e,x}\beta_2+d_2$, where
$\omega^a_{e,x}=\omega^1_{e,x}-\omega^2_{e,x}-\omega^3_{e,x}+\omega^4_{e,x}$,
$\omega^b_{e,x}=\omega^{+}_{e,x}-\omega^{-}_{e,x}+\omega^{\delta+}_{e,x}$,
$d_1=-\omega^{\delta-}_{e,x}\beta_1+
\omega^{\delta+}_{e,x}(\lambda_1-\beta_1)$ and
$d_2=-\omega^{\delta-}_{e,x}\lambda_2$.
Since $\lambda_1<\beta_1$ (see Proposition \ref{observt}), then
$d_1,d_2\leq 0$.
\item[(b)] if $\frac{\pi}{4}\leq\theta'<\theta<\theta''$,
\begin{eqnarray*}
\Delta'&=&(\omega^1_{e,x}-\omega^2_{e,x}-\omega^3_{e,x}+\omega^4_{e,x})\lambda_1+
(\omega^{+}_{e,x}-\omega^{-}_{e,x}+\omega^{\delta+}_{e,x})
\beta_1-\omega^{\delta-}_{e,x}\beta_1;\\
\Delta''&=&(-\omega^1_{e,x}+\omega^2_{e,x}+\omega^3_{e,x}-\omega^4_{e,x})\lambda_2+
(-\omega^{+}_{e,x}+\omega^{-}_{e,x}-\omega^{\delta+}_{e,x})
\beta_2-\omega^{\delta-}_{e,x}\lambda_2
+\omega^{\delta+}_{e,x}(\beta_2-\lambda_2).
\end{eqnarray*}
Hence, $\Delta'=\omega^a_{e,x}\lambda_1+\omega^b_{e,x}\beta_1+d'_1$
and $\Delta''=-\omega^a_{e,x}\lambda_2-\omega^b_{e,x}\beta_2+d'_2$,
where $d'_1=-\omega^{\delta-}_{e,x}\beta_1$ and
$d'_2=-\omega^{\delta-}_{e,x}\lambda_2+\omega^{\delta+}_{e,x}(\beta_2-\lambda_2)$.
By Proposition \ref{observt}, $\lambda_2>\beta_2$ and therefore,
$d'_1,d'_2\leq 0$.
\item[(c)] If $\theta'<\theta=\frac{\pi}{4}<\theta''$,
using similar arguments to the two previous cases, we obtain that
$\Delta'=\omega^a_{e,x}\lambda_1+\omega^b_{e,x}\beta_1+d_1$ and
$\Delta''=-\omega^a_{e,x}\lambda_2-\omega^b_{e,x}\beta_2+d'_2$. By
Proposition \ref{observt}, $\lambda_2>\beta_2$, $\beta_1>\lambda_1$
and therefore, $d_1,d'_2\leq 0$.\fin
\end{itemize}

\begin{teorema}
Let $e$, $x$ be the extreme points of a segment of length $\ell$
such that $x_1\neq a_{i1}$, $e_2\neq a_{i2}$, $\forall a_i\in A$
{ as well as $x$ and $e$ do not satisfy condition ii) of Theorem \ref{t21}}.
Then, $f(e',x')\leq f(e,x)$ or $f(e'',x'')\leq f(e,x)$, where $e'$,
$x'$, $e''$, $x''$ were given by (\ref{desplaza_2}).
\end{teorema}
\dem By Lemma \ref{prop42_2}, $f(x',e')-f(x,e)=\Delta'$ and
$f(x'',e'')-f(x,e)=\Delta'',$ where
$\Delta'=\omega^a_{e,x}\lambda_1+\omega^b_{e,x}\beta_1+d'$ and $
\Delta''=-\omega^a_{e,x}\lambda_2-\omega^b_{e,x}\beta_2+d''$, with
$d', d''\leq 0$.
\begin{enumerate}
\item If $\omega^b_{e,x}\leq 0$, we choose $\lambda_1=\lambda_2=\lambda$.
In this case, $\beta_2<\beta_1$ (Proposition \ref{observt}), then
$\Delta'+\Delta''=\omega^b_{e,x}(\beta_1-\beta_2)+d'+d''\leq 0$ and
the result follows.
\item If $\omega^b_{e,x}> 0$, we distinguish:
\begin{enumerate}
\item If $\omega^a_{e,x}\geq 0$, we choose $\beta_1=\beta_2=\beta$.
In this case, $\lambda_2>\lambda_1$ (Proposition \ref{observt}), then
$\Delta'+\Delta''=\omega^a_{e,x}(\lambda_1-\lambda_2)+d''+d''\leq
0$.
{
\item If $\omega^a_{e,x}< 0$, let $\theta$, $\theta'$ and $\theta''$ be the angles of the segments defined by
$e$ and $x$, $e'$ and $x'$, and $e''$ and $x''$, and the horizontal
line, respectively.
\begin{itemize}
\item[i)]If $\frac{\omega^b_{e,x}}{-\omega^a_{e,x}}< \tan\left(\theta\right)$, since $\theta-\theta'$ is small enough  then
$\frac{\omega^b_{e,x}}{-\omega^a_{e,x}} \le \tan\left(\frac{\theta+\theta'}{2}\right)$, equivalently,
$\frac{-\omega^a_{e,x}}{\omega^b_{e,x}} \ge
\frac{1}{\tan\left(\frac{\theta+\theta'}{2}\right)}$.
By Proposition \ref{observt}, we get that
$\frac{\beta_1}{\lambda_1}=\frac{1}{\tan\left(\frac{\theta+\theta'}{2}\right)}$.
Then $\frac{\beta_1}{\lambda_1} \le
\frac{-\omega^a_{e,x}}{\omega^b_{e,x}}$, i.e.,
$\omega^a_{e,x}\lambda_1+\omega^b_{e,x}\beta_1\le  0$. Therefore,
$\Delta'\leq 0$.
\item[ii)]If $\frac{\omega^b_{e,x}}{-\omega^a_{e,x}}> \tan\left(\theta\right)$, since $\theta''-\theta$ is small enough  then
$\frac{\omega^b_{e,x}}{-\omega^a_{e,x}}\geq
\tan\left(\frac{\theta+\theta''}{2}\right)$. Using similar arguments as in the previous case, we
obtain that $\Delta''\leq 0$.
\end{itemize}
}
\end{enumerate}
\end{enumerate}
Thus, we have found a movement where the objective function does not increase except for the case $\frac{\omega^b_{e,x}}{-\omega^a_{e,x}}=\tan\left(\theta\right)$
with $\omega^b_{e,x}>0$ and $\omega^a_{e,x}<0$, i.e., the case described by condition ii) of Theorem \ref{t21} or whenever one the endpoints of the optimal segment is an intersection point. Therefore the result follows.
 \fin




\end{document}